\documentclass[notitlepage,leqno,12pt]{article}

\usepackage{amsfonts}
\usepackage{amsmath}
\usepackage{amssymb}
\usepackage{color}

\definecolor{dark green}{rgb}{0,0.6,0}

\newcommand{\Chi}{\raise 1.5pt\hbox{$\chi$}}

\newcommand{\beq}{\begin{equation}}
\newcommand{\eeq}[1]{\label{#1}\end{equation}}

\newcommand{\avg}{\raise 0pt\hbox{$-$}\hskip -10.7pt\int}

\newcommand{\restrict}[1]{|_{\raise-2pt\hbox{$\scriptstyle #1$}}}

\newtheorem{lem}         {Lemma}[section]
\newtheorem{pro}    [lem]{Proposition}
\newtheorem{thm}    [lem]{Theorem}

\title{Existence of Solution for an Elliptic Problem with a Sublinear Term}
\date{}

\begin{document}

\small

\maketitle

\begin{center}
{\small
Anderson Araujo  \\ Universidade Federal de Vi\c{c}osa, Departamento de Matem\'atica, \\
Avenida Peter Henry Rolfs, \\ Vi\c{c}osa, MG, Brazil, CEP 36570-000 \\
E-mail: {\tt anderson.araujo@ufv.br} \\
\thanks{NOTE: This author was partially sponsored by FAPESP, Brazil, grant 2013/22328-8.}}
\end{center}

\begin{center}
{\small
Rafael Abreu  \\ Universidade Federal de Santa Catarina, Campus Blumenau, \\
Rua Pomerode, 710, \\ Blumenau, SC, Brazil, CEP 89065-300 \\
E-mail: {\tt rafael.abreu@ufsc.br}}
\end{center}

\noindent {\sc abstract}. In this work we prove the existence of a classical positive solution for an elliptic equation with a sublinear term. We use Galerkin approximations to show existence of such solution on bounded domains in $\mathbb{R}^N$.

\vspace{0.5cm}


\vspace{0.5cm}

\noindent {\sc keywords}. Elliptic problem, Galerkin method, bounded domain.

\vspace{0.5cm}

\section{Introduction}

In this paper, we study existence of solution for the problem
\begin{equation}\label{P1}
	\left\{
	\begin{array}{lcc}
	-\Delta v = \lambda v^{q} + f(v), &\textup{in}& \Omega, \\
	v>0 &\textup{in}& \Omega, \\
	v=0 &\textup{on}& \partial\Omega,
		\end{array}
		\right.
\end{equation}
where $\Omega \subset \mathbb{R}^N$, $N \geq 2$, is a bounded domain with smooth boundary, $\lambda>0$ is a parameter, $0<q<1$ and $f: \mathbb{R} \rightarrow \mathbb{R}$ is a continuous function satisfying
\begin{equation}\label{I2}
	0\leq f(s)s\leq C|s|^{p+1},
\end{equation}
where $1< p \leq \frac{N+2}{N-2}$ if $N\geq 3$ or $1< p$ if $N=2$.

Our main result in this paper is the following:
\begin{thm}\label{main result}
Suppose that $f:\mathbb{R} \to \mathbb{R}$ is a continuous function satisfying (\ref{I2}). Then, there exists $\lambda^*>0$ such that for every $\lambda \in (0,\lambda^*)$ the problem (\ref{P1}) has a positive solution $u \in C^{2,\gamma}(\overline{\Omega})$, for some $\gamma \in (0,1)$.
\end{thm}

Elliptic problems of the type
\begin{equation}\label{Pgeral}
	\left\{
	\begin{array}{lcc}
	-\Delta v = g(x,v) &\textup{in}& \Omega, \\
		v=0 &\textup{on}& \partial\Omega,
		\end{array}
		\right.
\end{equation}
where $g(x,v)$ is continuous and behaves like $v^{q} + v^p$ as $|v| \rightarrow +\infty$ have been extensively studied; see for example \cite{ambrosetti01,ambrosetti02,ABC} for a survey. One of the main results with nonlinearity combined effects of concave and convex was introduced in \cite{ABC}, namely, $g(x,u)=\lambda u^{q} + u^p$ with $0<q<1<p$.

We say that $g$ has sublinear growth at $+\infty$ if for every $\sigma \geq 0$ we have
\[
\lim_{|s| \rightarrow +\infty}\frac{g(x,s)}{|s|^{\sigma +1}}=0 \,\,\,\mbox{uniformily in} \,\, x
\]
and say that $g$ has superlinear growth at $+\infty$ if for every $\sigma\geq 0$ we have
\[
\lim_{|s| \rightarrow +\infty}\frac{g(x,s)}{|s|^{\sigma +1}}= \infty \,\,\,\mbox{uniformily in} \,\, x.
\]
We would like to righlight that the only assumptions which we assume are that $0<q<1$ and that $f$ is continuous and satisfies the growth condition (\ref{I2}). This way, the nonlinearity $g(s)=\lambda s^q + f(s)$ of problem (\ref{P1}) can have sublinear or superlinear growth at $+\infty$.

Most papers treat problem (\ref{Pgeral}) by means of variational methods, then it is usually assumed that $g$ has sublinear or superlinear growth and, sometimes, $sg(s) \geq c |s|^p$, where $c>0$ is a constant and $p>2$; see for example \cite{WZ}. Another common assumption on $g$ is the so-called Ambrosetti-Rabinowitz condition that means the following:
$$
\exists R>0 \mbox{ and } \theta > 2 \mbox{ such that } 0< \theta G(x,s) \leq s g(x,s) \,\, \forall |s| \geq R \mbox{ and } x \in \Omega,
$$
where $G(x,s)=\int_0^s g(x,\tau)d\tau$.

Even when the Ambrosetti-Rabinowitz condition can be dropped, it has to be assumed some condition to give compactness of Palais-Smale sequences or Cerami sequences. See for instance \cite{lg}, where they assume
$$
g: \overline \Omega \times \mathbb{R} \mbox{ is continuous and } g(x,0) = 0;
$$

$$
\exists t_0 > 0 \mbox{ and } M > 0 \mbox{ such that }
0 < G(x,s) \leq M g(x,s) \,\, \forall |s| \geq t_0 \mbox{ and } x \in \Omega;
$$

$$
0 < 2 G(x,s) \leq s g(x,s) \,\, \forall |s| \geq 0 \mbox{ and } x \in \Omega.
$$
See also \cite{MS}.

We are able to solve (\ref{P1}) under weaker assumptions by using the Galerkin method. For that matter we approximate $f$ by Lipschitz functions in Section \ref{s.2}. In Section \ref{s.3} we solve approximate problems. In Section \ref{s.4} we prove a regularity result to approximate problems. Section \ref{s.5} is devoted to prove Theorem \ref{main result}; in doing so we show that solutions $v_n$ of approximate problems are bounded away from zero and converge to a positive solution of (\ref{P1}).

At last in this introduction, we would like to emphasize that a similar approach was already used in \cite{af1}, but different to that, we do not assume that the nonlinearity $f$ is Lipschitz continuous.



	\section{Approximating functions}\label{s.2}

	In order to proof Theorem \ref{main result}, we make use of the following approximation result by Lipschitz functions, proved by Strauss in \cite{Strauss}.

	\begin{lem}\label{lemma1}
	Let $f:\mathbb{R} \to \mathbb{R}$ be a continuous function such that $sf(s)\geq 0$ for all $s \in \mathbb{R}$. Then, there exists a sequence $f_k:\mathbb{R} \to \mathbb{R}$ of continuous functions satisfying $sf_k(s)\geq 0$ and 
	\begin{description}
		\item[(i)] $\forall \, k \in \mathbb{N}$, $\exists c_k > 0$ such that $|f_k(\xi) - f_k(\eta)|\leq c_k|\xi - \eta|$,  for all $\xi, \eta \in \mathbb{R}$.
		\item[(ii)]  $(f_k)$ converges uniformly to $f$ in bounded subsets of $\mathbb{R}$.
	\end{description}
\end{lem}

The proof consists in considering the following family of approximation functions $f_k:\mathbb{R} \to \mathbb{R}$ defined by
\begin{eqnarray}\label{eq1}
f_k(s)=\displaystyle\left\{
	\begin{array}{lcc}
-k[G(-k-\frac{1}{k}) - G(-k)], &\textup{if}& s\leq -k,\\	
-k[G(s-\frac{1}{k}) - G(s)], &\textup{if}& -k\leq s \leq -\frac{1}{k},\\
k^2s[G(-\frac{2}{k}) - G(-\frac{1}{k})], &\textup{if}& -\frac{1}{k}\leq s\leq 0,\\
k^2s[G(\frac{2}{k}) - G(\frac{1}{k})], &\textup{if}& 0\leq s\leq \frac{1}{k},\\	
k[G(s+\frac{1}{k}) - G(s)], &\textup{if}& \frac{1}{k}\leq s \leq k,\\
k[G(k+\frac{1}{k}) - G(k)], &\textup{if}& s\geq k.\\
	\end{array}
		\right.
\end{eqnarray}
where $G(s)=\int_0^sf(\tau)d\tau$.

The sequence $(f_k)$ of the previous lemma has some additional properties.

\begin{lem}\label{lemma2}
Let $f: \mathbb{R} \to \mathbb{R}$ be a continuous function such that $sf(s)\geq 0$ for all $s \in \mathbb{R}$. Let us suppose that there exist constants $C>0$ and $1< p \leq \frac{N+2}{N-2}$ such that
\begin{equation}\label{eq2}
sf(s)\leq C|s|^{p+1} \ \ , \ \ \forall s \in \mathbb{R}. 	
\end{equation}
Then, the sequence $(f_k)_{k \in \mathbb{N}}$ from Lemma \ref{lemma1} satisfies

\begin{description}
	\item[(i)] $0\leq sf_k(s) \leq C_1|s|^{p+1}$ for all $|s|\geq \frac{1}{k}$,
	\item[(ii)]  $0\leq sf_k(s) \leq C_2|s|^2$ for all $|s|\leq \frac{1}{k}$,
\end{description}
where $C_1, C_2$ do not depend on $k$.
\end{lem}
\textbf{Proof:} Everywhere in this proof, the constant $C$ is the one given by (\ref{I2}).

\textit{\underline{First step}}: Suppose $-k\leq s \leq -\frac{1}{k}$.

By the mean value theorem, there exists $\eta\in (s-\frac{1}{k},s)$ such that
\[f_k(s)=-k[G(s-\frac{1}{k}) - G(s)]=-kG'(\eta)(s-\frac{1}{k}-s)=f(\eta)\]
and
\[sf_k(s)=sf(\eta).\]
As $s-\frac{1}{k}<\eta<s<0$ and $f(\eta)<0$, we have $sf(\eta) \leq \eta f(\eta)$. Therefore,
\[sf_k(s) \leq \eta f(\eta) \leq C|\eta|^{p+1} \leq C|s-\frac{1}{k}|^{p+1}\leq C(|s| +\frac{1}{k})^{p+1}\leq C2^{p+1}|s|^{p+1}.\]

\textit{\underline{Second step}}: Suppose $\frac{1}{k}\leq s \leq k$.

By the mean value theorem, there exist $\eta\in (s,s+\frac{1}{k})$ such that
\[f_k(s)=k[G(s+\frac{1}{k}) - G(s)]=kG'(\eta)(s+\frac{1}{k}-s)=f(\eta)\]
and
\[sf_k(s)=sf(\eta).\]
As $0<s <\eta <s+\frac{1}{k}$ and $f(\eta)>0$, we have $sf(\eta) \leq \eta f(\eta)$. Therefore,
\[sf_k(s) \leq \eta f(\eta) \leq C|\eta|^{p+1} \leq C|s+\frac{1}{k}|^{p+1} = C(|s| +\frac{1}{k})^{p+1}\leq C2^{p+1}|s|^{p+1}.\]

\textit{\underline{Third step}}: Suppose $|s|\geq k$.

Define 
\begin{eqnarray*}\label{eq3}
	f_k(s)=\displaystyle\left\{
	\begin{array}{lcc}
	-k[G(-k-\frac{1}{k}) - G(-k)], &\textup{if}& s\leq -k,\\
	k[G(k+\frac{1}{k}) - G(k)], &\textup{if}& s\geq k.\\
		\end{array}
		\right.
\end{eqnarray*}
 
If $s\leq -k$, by the mean value theorem, there exist $\eta\in (-k-\frac{1}{k},-k)$ such that
\[f_k(s)=k[G(-k-\frac{1}{k}) - G(-k)]=-kG'(\eta)(-k-\frac{1}{k}-(-k))=f(\eta)\]
and
\[sf_k(s)=sf(\eta).\]
As $-k-\frac{1}{k} <\eta <-k<0$ and $k<|\eta| < k + \frac{1}{k}$, we have $sf(\eta) = \frac{s}{\eta}\eta f(\eta)$. Therefore,
\[sf_k(s)=\frac{s}{\eta}\eta f(\eta) \leq \frac{|s|}{|\eta|}C|\eta|^{p+1} =\]
\[=C|s||\eta|^{p} \leq C|s|(k + \frac{1}{k})^{p} \leq C|s|(|s| + \frac{1}{k})^{p}\leq C2^{p}|s|^{p+1}.\]

If $s\geq k$, by the mean value theorem, there exist $\eta\in (k,k + \frac{1}{k})$ such that
\[f_k(s)=k[G(k+\frac{1}{k}) - G(k)]=kG'(\eta)(k+\frac{1}{k}-k)=f(\eta)\]
and
\[sf_k(s)=sf(\eta)=\frac{s}{\eta}\eta f(\eta) \leq \frac{|s|}{|\eta|}C|\eta|^{p+1} =\]
\[=C|s||\eta|^{p} \leq C|s|(k + \frac{1}{k})^{p} \leq C|s|(|s| + \frac{1}{k})^{p}\leq C2^{p}|s|^{p+1}.\]

\textit{\underline{Fourth step}}: Suppose $-\frac{1}{k}\leq s\leq \frac{1}{k}$.

Define 
\begin{eqnarray*}\label{eq3}
	f_k(s)=\displaystyle\left\{
	\begin{array}{lcc}
	k^2s[G(-\frac{2}{k}) - G(-\frac{1}{k})], &\textup{if}& -\frac{1}{k}\leq s\leq 0,\\
	k^2s[G(\frac{2}{k}) - G(\frac{1}{k})], &\textup{if}& 0\leq s\leq \frac{1}{k}.\\
		\end{array}
		\right.
\end{eqnarray*}
 
If $-\frac{1}{k}\leq s\leq 0$, by the mean value theorem, there exists $\eta\in (-\frac{2}{k},-\frac{1}{k})$ such that
\[f_k(s)=k^2s[G(-\frac{2}{k}) - G(-\frac{1}{k})]=k^2sG'(\eta)(-\frac{2}{k}-(-\frac{1}{k}))=-ksf(\eta).\]
Therefore,
\[sf_k(s)=-ks^2f(\eta)=-k\frac{s^2}{\eta}\eta f(\eta) \leq k\frac{s^2}{|\eta|}\eta f(\eta)\]
\[\leq Ck|s|^2|\eta|^{p} \leq Ck|s|^2(\frac{2}{k})^{p} \leq C2^{p}|s|^2.\]

If $0\leq s \leq \frac{1}{k}$, by the mean value theorem, there exist $\eta\in (\frac{1}{k},\frac{2}{k})$ such that
\[f_k(s)=k^2s[G(\frac{2}{k}) - G(\frac{1}{k})]=k^2sG'(\eta)(\frac{2}{k}-\frac{1}{k})=ksf(\eta).\]
Therefore,
\[sf_k(s)=ks^2f(\eta)=k\frac{s^2}{|\eta|}\eta f(\eta) \leq \]
\[\leq Ck|s|^2|\eta|^{p} \leq Ck|s|^2(\frac{2}{k})^{p} \leq C2^{p}|s|^2.\]
The proof of the lemma follows by taking $C_1=C2^{p+1}$ and $C_2=C2^{p}$, where $C$ is like in (\ref{eq2}).

\section{Approximate problem} \label{s.3}

In order to prove Theorem \ref{main result}, we first study the auxiliary problem
\begin{equation}\label{P2}
	\left\{
	\begin{array}{lcc}
	-\Delta v = \lambda v^{q} + f_n(v) + \frac{1}{n} &\textup{in}& \Omega,\\
	v>0 &\textup{in}& \Omega,\\
	v=0 &\textup{on}& \partial\Omega,
		\end{array}
		\right.
\end{equation}
where $0<q<1$, $\lambda > 0$ is a parameter and $f_n: \mathbb{R} \rightarrow \mathbb{R}$ is a function of the sequence given by Lemma \ref{lemma1} and Lemma \ref{lemma2}.


We will use the Galerkin method together with the following fixed point theorem, see \cite{Strauss} and \cite[Theorem 5.2.5]{k}. A similar approach was already used in \cite{af1}.

\begin{pro}\label{prop1}
Let $F: \mathbb{R}^d \rightarrow \mathbb{R}^d$ be a continuous function such that $\left\langle F(\xi),\xi\right\rangle\geq 0$ for every $\xi \in \mathbb{R}^d$ with $|\xi|=r$ for some $r>0$. Then, there exists $z_0$ in the closed ball $\overline{B}_r(0)$ such that $F(z_0)=0$. 
\end{pro}

The main result in this section is the following theorem.
\begin{thm}\label{teo aux}
There exists $\lambda^*>0$ and $n^* \in \mathbb{N}$ such that (\ref{P2}) has a weak positive solution for all $\lambda \in (0,\lambda^*)$ and $n\geq n^*$.

\end{thm}
\textbf{Proof:} Fix $\mathcal{B}=\{w_1,w_2,\dots,w_m,\dots\}$ a orthonormal basis of $H_0^1(\Omega)$ and define
\[W_m=[w_1,w_2,\dots,w_m],\]
to be the space generated by $\{w_1,w_2,\dots,w_m\}$. Define the function $F:\mathbb{R}^m \to \mathbb{R}^m$ such that $F(\xi)=(F_1(\xi),F_2(\xi),\dots, F_m(\xi))$, where
\[F_j(\xi)=\int_{\Omega}\nabla v\nabla w_j - \lambda\int_{\Omega}(v_+)^{q}w_j - \int_{\Omega}f_n(v_+)w_j - \frac{1}{n}\int_{\Omega}w_j, \,\,\,j=1,2,\dots,m\]
and let $v=\sum_{i=1}^m\xi_iw_i$. Therefore, 
\begin{equation}\label{eq4}
\left\langle F(\xi),\xi \right\rangle=\int_{\Omega}|\nabla v|^2 - \lambda\int_{\Omega}(v_+)^{q+1} - \int_{\Omega}f_n(v_+)v_+ - \frac{1}{n}\int_{\Omega}v.	
\end{equation}

Given $v \in W_m$ we define 
\[\Omega^+_n=\{x \in \Omega : |v(x)|\geq \frac{1}{n}\}\] 
and 
\[\Omega^-_n=\{x \in \Omega : |v(x)|< \frac{1}{n}\}.\] 
Thus we rewrite (\ref{eq4}) as
\[\left\langle F(\xi),\xi\right\rangle = \left\langle F(\xi),\xi\right\rangle_P + \left\langle F(\xi),\xi\right\rangle_{N},\]
where
\[\left\langle F(\xi),\xi\right\rangle_P=\int_{\Omega^+_n}|\nabla v|^2 - \lambda\int_{\Omega^+_n}(v_+)^{q+1} - \int_{\Omega^+_n}f_n(v_+)v_+ - \frac{1}{n}\int_{\Omega^+_n}v	\]
and
\[\left\langle F(\xi),\xi\right\rangle_{N}=\int_{\Omega^-_n}|\nabla v|^2 - \lambda\int_{\Omega^-_n}(v_+)^{q+1} - \int_{\Omega^-_n}f_n(v_+)v_+ - \frac{1}{n}\int_{\Omega^-_n}v.\]

\textit{\underline{Step 1}}. Since $0<q<1$, then
\begin{equation}\label{eq4.1}
	\int_{\Omega^+_n}(v_+)^{q+1} \leq \int_{\Omega}|v|^{q+1} = \|v\|^{q+1}_{L^{q+1}(\Omega)}\leq C_1\|v\|^{q+1}_{H^1_0(\Omega)}.
\end{equation}
By virtue of ($i$) Lemma \ref{lemma2} we get
\begin{equation}\label{eq4.2}
\begin{array}{rcl}
\displaystyle\int_{\Omega^+_n}f_n(v_+)v_+ &\leq&\displaystyle C\int_{\Omega}|v_+|^{p+1}dx\leq C_2\|v\|^{p+1}_{H^1_0(\Omega)}.
\end{array}
\end{equation}

It follows from (\ref{eq4.1}) and (\ref{eq4.2}) that
\begin{equation}\label{eq5}
\begin{array}{rcl}
	\left\langle F(\xi),\xi\right\rangle_P &\geq &\displaystyle \int_{\Omega^+_n}|\nabla v|^2 - \lambda C_1\|v\|^{q+1}_{H^1_0(\Omega)}\\
	&-& \displaystyle C_2\|v\|^{p+1}_{H^1_0(\Omega)} - \frac{C_3}{n}\|v\|_{H^1_0(\Omega)},
	\end{array}
\end{equation}
where $C_1$, $C_2$ and $C_3$ depends on $C$ and $|\Omega|$.

\textit{\underline{Step 2}}. Since $0<q<1$, then 
\begin{equation}\label{eq10}
\int_{\Omega^-_n}(v_+)^{q+1} \leq \int_{\Omega^-_n}|v|^{q+1} \leq |\Omega|\frac{1}{n^{q+1}}.
\end{equation}
By virtue of ($ii$) Lemma \ref{lemma2} we get
\begin{equation}\label{eq11}
\int_{\Omega^-_n}f_n(v_+)v_+ \leq C\int_{\Omega^-_n}|v_+|^2dx\leq C|\Omega|\frac{1}{n^2}.
\end{equation}
It follows from (\ref{eq10}) and (\ref{eq11}) that
\begin{equation}\label{eq12}
	\left\langle F(\xi),\xi\right\rangle_N \geq \int_{\Omega^-_n}|\nabla v|^2 - \lambda |\Omega|\frac{1}{n^{q+1}} - C|\Omega|\frac{1}{n^2} - |\Omega|\frac{1}{n^2}.
\end{equation}

It follows from (\ref{eq5}) and (\ref{eq12}) that
\begin{eqnarray*}\label{eq13}
\begin{array}{rcl}
\displaystyle\left\langle F(\xi),\xi\right\rangle  &\geq& \displaystyle \|v\|^2_{H^1_0(\Omega)} - \lambda C_1\|v\|^{q+1}_{H^1_0(\Omega)} - C_2\|v\|^{p+1}_{H^1_0(\Omega)} \\
& - &\displaystyle \frac{C_3}{n}\|v\|_{H^1_0(\Omega)} - \lambda |\Omega|\frac{1}{n^{q+1}} - C|\Omega|\frac{1}{n^2} - |\Omega|\frac{1}{n^2}.
\end{array}
\end{eqnarray*}
Assume now that $\|v\|_{H^1_0(\Omega)}=r$ for some $r>0$ to be fixed later. Hence,
\[\left\langle F(\xi),\xi\right\rangle \geq r^2 - \lambda C_1r^{q+1} - C_2r^{p+1} - \frac{C_3}{n}r - \lambda |\Omega|\frac{1}{n^{q+1}} - C|\Omega|\frac{1}{n^2} - |\Omega|\frac{1}{n^2}.\]
We want to choose $r$ such that
\[r^2 - C_2r^{p+1}\geq \frac{r^2}{2},\]
in other words,
\[r\leq \frac{1}{(2C_2)^{\frac{1}{p-1}}}.\]
Choosing $r=\frac{1}{2(2C_2)^{\frac{1}{p-1}}}$, we obtain
\[\left\langle F(\xi),\xi\right\rangle \geq \frac{r^2}{2} - \lambda C_1r^{q+1} - \frac{C_3}{n}r - \lambda |\Omega|\frac{1}{n^{q+1}} - C|\Omega|\frac{1}{n^2} - |\Omega|\frac{1}{n^2}.\]
Now, defining  $\rho=\frac{r^2}{2} - \lambda C_1r^{q+1}$, we choose $\lambda^*>0$ such that $\rho>0$ for $\lambda < \lambda^*$. 
Therefore, we choose $\lambda^*=\frac{r^{1-q}}{4C_1}$. Now we choose $n^* \in \mathbb{N}$ such that
\[\frac{C_3}{n}r + \lambda |\Omega|\frac{1}{n^{q+1}} + C|\Omega|\frac{1}{n^2} + |\Omega|\frac{1}{n^2} <\frac{\rho}{2},\]
for every $n\geq n^*$. Let $\xi \in \mathbb{R}^m$, such that $|\xi|=r$, then for $\lambda < \lambda^*$ and  $n\geq n^*$ we obtain
\begin{eqnarray*}\label{eq8}
	\left\langle F(\xi),\xi\right\rangle \geq \frac{\rho}{2}>0.
\end{eqnarray*}

Since $f_n$ is a Lipschitz continuous function for every $n$, by standard arguments it is shown that $F$ is continuous, that is, give $(x_k)$ in $\mathbb{R}^m$ and $x \in \mathbb{R}^m$ such that $x_k \rightarrow x$ we obtain $F(x_k) \rightarrow F(x)$.

Therefore, by Proposition \ref{prop1} for all $m \in \mathbb{N}$ there exists $y \in \mathbb{R}^m$ with $|y|\leq r$ such that $F(y)=0$, that is, there exists $v_m \in W_m$ verifying $\|v_m\|_{H^1_0(\Omega)}\leq r$, for every $m \in \mathbb{N}$ and such that
\begin{eqnarray*}\label{eq15}
\int_{\Omega}\nabla v_m\nabla w = \lambda\int_{\Omega}(v_{m+})^{q}w + \int_{\Omega}f_n(v_{m+})w + \frac{1}{n}\int_{\Omega}w, \,\,\, \forall \, w \in W_m.
\end{eqnarray*}

Since $W_m \subset H^1_0(\Omega)$, $\forall \, m \in \mathbb{N}$,  and  $r$ does not depend on $m$, then $(v_m)$ is a bounded sequence of $H^1_0(\Omega)$. Then, for some subsequence, there exists $v=v_n \in H^1_0(\Omega)$ such that
 \begin{equation}\label{eq16}
v_m \rightharpoonup v \,\,\, \mbox{weakly in} \,\,\, H^1_0(\Omega)
\end{equation}
and
 \begin{equation}\label{eq16.1}
v_m \to v \,\,\, \mbox{in} \,\,\, L^2(\Omega)\,\,\, \mbox{and a.e. in}\,\,\,\Omega.
\end{equation}
Fixing	$k \in \mathbb{N}$ and for every $m$ such that $m\geq k$ we obtain
\begin{equation}\label{eq17}
\int_{\Omega}\nabla v_m\nabla w_k = \lambda\int_{\Omega}(v_{m+})^{q}w_k + \int_{\Omega}f_n(v_{m+})w_k + \frac{1}{n}\int_{\Omega}w_k, \,\,\, \forall \, w_k \in W_k.
\end{equation}

Now, as $g:H^1_0(\Omega) \rightarrow \mathbb{R}$ defined by $g(u)=\int_{\Omega}\nabla u\nabla w_k$, for every $u \in H^1_0(\Omega)$, we have that $g$ is a continuous linear functional. It follows from (\ref{eq16}) that
\begin{equation}\label{eq18}
\int_{\Omega}\nabla v_m\nabla w_k \rightarrow \int_{\Omega}\nabla v\nabla w_k \,\,\,\mbox{as}\,\,\, m \to \infty
\end{equation}
and by (\ref{eq16.1}), we obtain
\begin{equation}\label{eq18.1}
\int_{\Omega}f_n(v_{m+})w_k \to \int_{\Omega}f_n(v_{+})w_k \,\,\,\mbox{as}\,\,\, m \to \infty.
\end{equation}
Indeed, by Lemma \ref{lemma1} ($ii$) it follows that $|f_n(v_{m+}) - f_n(v_{+})|\leq c_n|v_{m+}-v_{+}|$, hence
\[\left|\int_{\Omega}f_n(v_{m+})w_k - \int_{\Omega}f_n(v_{+})w_k\right| \leq c_n\|w_k\|_{L^2(\Omega)}\|v_{m}-v\|_{L^2(\Omega)} \,\,\,\mbox{as}\,\,\, m \to \infty,\]
and then, (\ref{eq16.1}) implies (\ref{eq18.1}).
By (\ref{eq16}), (\ref{eq18.1}) and Sobolev compact embedding, letting $m \to \infty$, we obtain
\begin{equation}\label{eq19}
\lambda\int_{\Omega}(v_{m+})^{q}w_k + \int_{\Omega}f_n(v_{m+})w_k + \frac{1}{n}\int_{\Omega}w_k \rightarrow \lambda\int_{\Omega}(v_{+})^{q}w_k + \int_{\Omega}f_n(v_{+})w_k + \frac{1}{n}\int_{\Omega}w_k.
\end{equation}
By (\ref{eq17}), (\ref{eq18}), (\ref{eq19}) and by the uniqueness of the limit, we obtain
\begin{eqnarray*}\label{eq20}
\int_{\Omega}\nabla v\nabla w_k = \lambda\int_{\Omega}(v_{+})^{q}w_k + \int_{\Omega}f_n(v_{+})w_k + \frac{1}{n}\int_{\Omega}w_k, \,\,\, \forall \, w_k \in W_k.
\end{eqnarray*}
For density of $[W_k]_{k \in \mathbb{N}}$ in $H^1_0(\Omega)$ and by linearity, we conclude that
\begin{equation}\label{eq21}
\int_{\Omega}\nabla v\nabla w = \lambda\int_{\Omega}(v_{+})^{q}w + \int_{\Omega}f_n(v_{+})w + \frac{1}{n}\int_{\Omega}w, \,\,\, \forall \, w \in H^1_0(\Omega).
\end{equation}
Furthermore, $v\geq 0$ in $\Omega$. In fact, as $v_- \in H^1_0(\Omega)$, we obtain from (\ref{eq21}) that
\[\int_{\Omega}\nabla v\nabla v_- = \lambda\int_{\Omega}(v_{+})^{q}v_- + \int_{\Omega}f_n(v_{+})v_- + \frac{1}{n}\int_{\Omega}v_-.\]
Hence, we have from Lemma \ref{lemma1} that
\[0\geq - \|v_-\|^2_{H^1_0(\Omega)} = \int_{\Omega}\nabla v\nabla v_- = \int_{\Omega}f_n(v_{+})v_- + \frac{1}{n}\int_{\Omega}v_- \geq 0,\]
with the result that $\|v_-\|_{H^1_0(\Omega)}=0$, that is, $v_-(x)=0$ a.e. in $\Omega$. Therefore, $v(x)=v_+(x)\geq 0$ a.e. in $\Omega$ and we conclude the proof of the theorem.

\section{Regularity of Solution of the Approximate Problem }\label{s.4}

In this section, we show that all weak solutions of the problem (\ref{P2}) are regular. Let $v \in H^1_0(\Omega)$ be a weak solution of the problem (\ref{P2}) and define
\begin{eqnarray*}
g(x) := \lambda v^q(x) + f_n(v(x)) + \frac{1}{n}.
\end{eqnarray*}
We have that
\begin{eqnarray}\label{reg1}
|g| \leq \lambda |v|^q + |f_n(v)| + \frac{1}{n}.
\end{eqnarray}
Notice that
\begin{eqnarray}\label{reg2}
|v|^q \leq 1 + |v|^{t-1},
\end{eqnarray}
where $2 \leq t \leq 2^*$. Here, $2^*$ is the critical Sobolev exponent, that is,
\begin{eqnarray*}
2^* = \frac{2N}{N-2}.
\end{eqnarray*}
Furthermore, since $f_n : \mathbb{R} \rightarrow \mathbb{R}$ is a Lipschitz continuous function and $f_n(0) = 0$, we have for each $n \in \mathbb{N}$ that
\begin{eqnarray*}
|f_n(v)| \leq C_n |v|,
\end{eqnarray*}
and consequently,
\begin{eqnarray}\label{reg3}
|f_n(v)| \leq C_n (1 + |v|^{t-1}),
\end{eqnarray}
where $2 \leq t \leq 2^*$. This way, by combining (\ref{reg1}), (\ref{reg2}) and (\ref{reg3}), we obtain
\begin{eqnarray}\label{reg4}
|g| \leq C_1 + C_2 |v|^{t-1},
\end{eqnarray}
where
\begin{eqnarray*}
C_1:= \lambda + C_n + \frac{1}{n}
\end{eqnarray*}
and
\begin{eqnarray*}
C_2:= \lambda + C_n.
\end{eqnarray*}
Then, using (\ref{reg4}) and well-known Bootstrap arguments, similar to those found in \cite{kavian}, we conclude that $v \in C^{2,\gamma}(\overline{\Omega})$, for some $\gamma \in (0,1)$.

\section{Proof of the Theorem \ref{main result}}\label{s.5}

In this section, we demonstrate Theorem \ref{main result}. The following lemma of \cite[Theorem 1.1]{Strauss} is used to show that $v_n$ converges to a solution $v$ of (\ref{P1}).
\begin{lem}\label{Strauss2}
Let $\Omega$ be a bounded open set in $\mathbb{R}^N$, $u_k : \Omega \to \mathbb{R}$ be a sequence of functions and $g_k : \mathbb{R} \to \mathbb{R}$ be a sequence of functions such that $g_k(u_k)$ are measurable in $\Omega$ for every $k \in \mathbb{N}$. Assume that $g_k(u_k) \to v$ a.e. in $\Omega$ and $\int_{\Omega}|g_k(u_k)u_k|dx<C$ for a constant $C$ independent of $k$. Suppose that for every bounded set $B \subset \mathbb{R}$ there is a constant $C_B$ depending only on $B$ such that $|g_k(x)|\leq C_B$, for all $x \in B$ and $k \in \mathbb{N}$. Then $v \in L^1(\Omega)$ and $g_k(u_k) \to v$ in $L^1(\Omega)$.
\end{lem}

Since $v \in C^{2,\gamma}(\overline{\Omega})$, $\gamma \in (0,1)$, satisfies $v\geq 0$ and
\begin{eqnarray*}
- \Delta v = \lambda v^q + f_n(v) + \frac{1}{n},
\end{eqnarray*}
it follows by assumptions on $f_n$ that
\begin{eqnarray*}
-\Delta v \geq 0.
\end{eqnarray*} 
Then, by Maximum Principle, we have $v>0$ in $\Omega$, that is, $v$ is a solution of the problem (\ref{P2}). For each $n \in \mathbb{N}$, let us denote by $v_n$ the solution of (\ref{P2}). It follows from (\ref{eq16}) that 
\begin{eqnarray*}\label{eq22}
v_m^{(n)} \rightharpoonup v_n \,\,\, \mbox{weakly in} \,\,\, H^1_0(\Omega) \ \ \mbox{as} \ \ m \rightarrow \infty,
\end{eqnarray*}
where, for each $n \in \mathbb{N}$, $(v_m^{(n)})_{m \in \mathbb{N}}$ is a sequence in $H^1_0(\Omega)$ satisfying
\begin{eqnarray*}
||v_m^{(n)}|| \leq r, \ \ \forall m \in \mathbb{N}.
\end{eqnarray*}
Then,
\[\|v_n\| \leq \liminf_{m \rightarrow \infty}\|v_m^{(n)}\|\leq r, \,\,\, \forall\, n \in \mathbb{N}.\]
Since $r$ does not depend on $n$, there exists $v \in H^1_0(\Omega)$ such that
\begin{eqnarray*}\label{eq23}
v_n \rightharpoonup v \,\,\, \mbox{weakly in} \,\,\, H^1_0(\Omega).
\end{eqnarray*}
By compact embedding, up to a subsequence, we have
\[v_n \rightarrow v \,\,\mbox{in}\,\,L^s(\Omega),\,\,\mbox{for} \ \ 1\leq s <2^* \,\, \mbox{if}\,\,N\geq 3 \,\, \mbox{or for}\,\, 1\leq s < +\infty\,\, \mbox{if}\,\, N=2,\]
and then, up to a subsequence,
\begin{description}
	\item[i)] $v_n(x) \rightarrow v(x)$ a.e. in $\Omega$;
	\item[ii)] $|v_n(x)|\leq h(x)$, $\forall\, n \in \mathbb{N}$ a.e. in $\Omega$, for some $h \in L^s(\Omega)$.
\end{description}

Notice that the following inequality holds:
\begin{eqnarray*}\label{P3}
	\left\{
	\begin{array}{lcc}
	-\Delta v_n \geq \lambda v_n^{q} , &\textup{in}& \Omega,\\
	v_n>0 &\textup{in}& \Omega,\\
	v_n=0 &\textup{on}& \partial\Omega.
		\end{array}
		\right.
\end{eqnarray*}
This way, considering $w_n=\lambda^{\frac{1}{q-1}}v_n$, we obtain
\[-\Delta\left(\frac{w_n}{\lambda^{\frac{1}{q-1}}}\right)\geq \lambda\left(\frac{w_n}{\lambda^{\frac{1}{q-1}}}\right)^q,\]
and consequently,
\[-\Delta\,w_n\geq w_n^q.\]
Let us denote by $\widetilde{w}$ the unique solution of the problem
\[
	\left\{
	\begin{array}{lcc}
	-\Delta\widetilde{w} = \widetilde{w}^{q} , &\textup{in}& \Omega,\\
	\widetilde{w}>0 &\textup{in}& \Omega,\\
	\widetilde{w}=0 &\textup{on}& \partial\Omega.
		\end{array}
		\right.
\]
The existence and uniqueness of such solution is proved in \cite{bo}. By Lemma 3.3 of \cite{ABC}, it follows that $w_n \geq \widetilde{w}$, $\forall\, n \in \mathbb{N}$, that is, 
\begin{equation}\label{eq24}
v_n(x) \geq \lambda^{\frac{1}{1-q}}\widetilde{w}(x), \,\,\mbox{a.e. in}\,\,\Omega, \forall\, n \in \mathbb{N}.
\end{equation}
Taking the limit as $n \rightarrow +\infty$ in (\ref{eq24}), we obtain
\[v(x) \geq \lambda^{\frac{1}{1-q}}\widetilde{w}(x), \,\,\mbox{a.e. in}\,\,\Omega\]
and hence $v>0$ a.e. in $\Omega$.

Recall that, from (\ref{eq21}),
\begin{eqnarray*}\label{eq25}
\int_{\Omega}\nabla v_n\nabla w = \lambda\int_{\Omega}(v_n)^{q}w + \int_{\Omega}f_n(v_n)w + \frac{1}{n}\int_{\Omega}w, \,\,\, \forall \, w \in H^1_0(\Omega) ,
\end{eqnarray*}
and using that $v_n$ is a classical solution we have
\begin{equation}\label{eq26}
-\Delta\,v_n = \lambda(v_n)^{q} + f_n(v_n) + \frac{1}{n} \,\,\, \mbox{in}\,\, L^2(\Omega).
\end{equation}
Since
\[
v_n \rightarrow v\,\, \mbox{a.e. in}\,\, \Omega,
\]
we have
\begin{equation}\label{eq26.2}
f_n(v_n(x)) \rightarrow f(v(x))\,\, \mbox{a.e. in}\,\, \Omega
\end{equation}
by the uniform convergence of Lemma \ref{lemma1} ($ii$).

Multiplying the equation (\ref{eq26}) by $w=v_n$ and since $v_n$ is bounded in $H_0^1(\Omega)$ we obtain
\begin{equation}\label{eq27}
\int_{\Omega}f_n(v_n)v_ndx \leq C,
\end{equation}
for every $n \in \mathbb{N}$, where $C>0$ is a constant independent of $n$. By (\ref{eq26.2}), (\ref{eq27}) and by the expression  of $f_n$ defined in (\ref{eq1}), the assumptions of Lemma \ref{Strauss2} are satisfied implying
\begin{eqnarray*}\label{eq28}
f_n(v_n) \rightarrow f(v)\,\, \mbox{strongly in}\,\, L^1(\Omega).
\end{eqnarray*}
Multiplying (\ref{eq26}) by $w \in \mathcal{D}(\Omega)$, integrating on $\Omega$ and using the previous convergences, we have
\begin{equation}\label{eq29}
-\Delta\,v = \lambda\,v^{q} + f(v) \,\,\, \mbox{in}\,\, \mathcal{D}'(\Omega).
\end{equation}

Since $f(v) \in L^{\frac{p+1}{p}}(\Omega)$ and $\lambda\,v^{q} \in L^{\frac{p+1}{p}}(\Omega)$, we conclude from (\ref{eq29}) that $v \in H^1_0(\Omega)\cap W^{2,\frac{p+1}{p}}(\Omega)$ and 
\[-\Delta v= \lambda\,v^{q} + f(v)\]
in the strong sense. Notice that the assumption (\ref{I2}) implies that
\begin{eqnarray*}
|f(s)| \leq C|s|^{t-1},
\end{eqnarray*}
where $2 \leq t \leq 2^*$. Thus, using well-known Bootstrap arguments, we conclude that $v \in C^{2,\gamma}(\overline{\Omega})$, for some $\gamma \in (0,1)$, and it is a classical positive solution of problem (\ref{P1}).



\bibliographystyle{amsplain}

\begin{thebibliography}{10}

\bibitem{af1}
C. O. Alves and D. G. de Figueiredo, \emph{Nonvariational Elliptic Systems via Galerkin Methods}, D. Haroske, T. Runst and H. J. Schmeisser (eds.) Function Spaces, Differential Operators and Nonlinear Analysis. The Hans Triebel Anniversary Volume, 2003

\bibitem{ambrosetti01}
A. Ambrosetti, \emph{Critical Points and Nonlinear Variational Problems}, Bull. Soc. Math. France 120, Memoire No. 49 (1992).

\bibitem{ambrosetti02}
A. Ambrosetti and M. Badiale, \emph{The Dual Variational Principle and Elliptic Problems with Discontinuous Nonlinearities}, J. Math. Anal. Appl. 140 (1989), 363-373.

\bibitem{ABC} A. Ambrosetti, H. Brezis, G. Cerami, \textit{Combined Effects of Concave and Convex Nonlinearities in Some Elliptic Problems},  Journal of Functional Analysis , \textbf{122}, (1994), 519--543.

\bibitem{bo}
H. Brezis and L. Oswald, \emph{Remarks on sublinear elliptic equations}. Nonlinear Analysis
TMA. \textbf{10} 55--64 (1986)

\bibitem{lg}
N. Lam and G. Lu, \emph{Elliptic equations and systems with subcritical and critical exponential growth without the Ambrosetti-Rabinowitz condition}. J. Geom. Anal. \textbf{24} 118--143 (2014)

\bibitem{kavian} O. Kavian, \textit{Introction à la théorie de Points Critiques}, vol. 13, Springer-Verlag, 1993

\bibitem{k}
S. Kesavan, Topics in functional analysis and applications, John Wiley $\&$ Sons (1989).

\bibitem{MS}
O. H. Miyagaki and M. A. S. Souto, \textit{Superlinear problems without Ambrosetti and Rabinowitz growth condition}, J. Differential Equations 245 (2008) 3628-3638.




\bibitem{Strauss}
W. A. Strauss, \emph{On weak solutions of semilinear hyperbolic equations}, An. Acad. Brasil. Ci\^enc. \textbf{42} 645--651 (1970)


\bibitem{WZ} M. Willem and W. Zou, \emph{On a Schrödinger equation with periodic potential and spectrum point zero},
Indiana Univ. Math. J. 52 (1), (2003), 109-132. 

\end{thebibliography}

\end{document}